\documentstyle[12pt,thmsa,sw20lart]{article}
\input{tcilatex}
\begin{document}

\author{{\small TEODOR OPREA} \\
%EndAName
University of Bucharest,\\
Faculty of Mathematics and Informatics\\
$14$ Academiei St., code 010014, Bucharest, Romania\\
e-mail: teodoroprea@yahoo.com}
\title{{\large A NEW OBSTRUCTION TO MINIMAL ISOMETRIC IMMERSIONS INTO A REAL SPACE
FORM}}
\date{}
\maketitle

\begin{center}
\medskip {\bf Abstract}
\end{center}

In the theory of minimal submanifolds, the following problem is funda-

mental: {\sl when does a given Riemannian manifold admit (or does not ad-}

{\sl mit) a minimal isometric immersion into an Euclidean space of arbitrary}

{\sl dimension? }S.S. Chern, in his monograph {\sl Minimal submanifolds in a
Rie-}

{\sl mannian manifold}, remarked that the result of Takahashi ({\sl the
Ricci ten-}

{\sl sor of a minimal submanifold into a Euclidean space is negative semide-}

{\sl finite)} was the only known Riemannian obstruction to minimal isometric

immersions in Euclidean spaces. A second solution to this problem was

obtained by B.Y. Chen as an immediate application of his fundamental

inequality [1]: {\sl the scalar curvature and the sectional curvature of a
mini-}

{\sl mal submanifold into a Euclidean space} {\sl satisfies the inequality }$%
\tau \leq k.$ In

this paper we prove that the sectional curvature of a minimal submanifold

into a Euclidean space also satisfies the inequality $k\leq -\tau .$\newline
\medskip


\medskip

\medskip

{\it Key\ words: }constrained{\bf \ }maximum, Chen's inequality, minimal
submani-

\qquad \qquad \quad \thinspace \thinspace folds

\bigskip

\begin{center}
{\bf 1. OPTIMIZATIONS\ ON\ RIEMANNIAN\ MANIFOLDS}
\end{center}

~\quad \medskip

Let $(N,\widetilde{g})$ be a Riemannian manifold, $(M,g)$ a Riemannian
submanifold, and $f\in {\cal F}(N).$ To these ingredients we attach the
optimum problem\newline

\[
(1)\text{ }\min\limits_{x\in M}f(x). 
\]

Let's remember the result obtained in [6].

\ 

T{\footnotesize HEOREM} 1.1.{\sl \ If }$x_{0}\in M${\sl \ is the solution of
the problem }$(1)${\sl , then}

{\sl \ }

i){\sl \ }$($grad f$)(x_{0})${\sl \ }$\in T_{x_{0}}^{\perp }M,$

{\sl \ }

ii){\sl \ the bilinear form }

\[
\begin{array}{c}
\alpha {\sl \ :\ }T_{x_{0}}M\times T_{x_{0}}M{\sl \ }\rightarrow {\sl \ }R,
\\ 
\alpha (X,Y)=\text{Hess}_{f}(X,Y)+\widetilde{g\text{ }}(h(X,Y),(\text{grad f}%
)(x_{0}))
\end{array}
\]
{\sl is positive semidefinite, where }$h${\sl \ is the second fundamental
form of the submanifold }$M${\sl \ in }$N.$

\begin{center}
\bigskip

{\bf 2.} {\bf THE CHEN'S INEQUALITY}
\end{center}

\ 

Let $(M,g)$ be a Riemannian manifold of dimension $n$, and $x$ a point in $%
M. $ We consider the orthonormal frame $\{e_{1},e_{2},...,e_{n}\}$ in $%
T_{x}M.$

\medskip

The scalar curvature at $x$ is defined by

\[
\tau =\dsum\limits_{1\leq i<j\leq n}R(e_i,e_j,e_i,e_j). 
\]

\ 

We denote 
\[
\delta _{M}=\tau -\min (k), 
\]
where $k$ is the sectional curvature at the point $x.$ The invariant $\delta
_{M}$ is called the Chen's invariant of Riemannian manifold $(M,g).$

\ 

The Chen's invariant can be estimated, if $(M,g)$ is a Riemannian
submanifold in a real space form $\widetilde{M}(c)$, varying with $c$ and
the mean curvature of $M$ in $\widetilde{M}(c).$

\medskip \medskip

\bigskip T{\footnotesize HEOREM} 2.1. {\sl Consider ( }$\widetilde{M}(c),%
\widetilde{g}${\sl ) a real space form of dimension }$m${\sl , }$M\subset 
\widetilde{M}(c)${\sl \ a Riemannian submanifold of dimension }$n\geq 3${\sl %
. The Chen's invariant of }$M${\sl \ satisfies} 
\[
\delta _{M}\leq \frac{n-2}{2}\{\frac{n^{2}}{n-1}\left\| H\right\|
^{2}+(n+1)c\}, 
\]
{\sl where }$H${\sl \ is the mean curvature vector of submanifold }$M${\sl \
in }$\widetilde{M}(c)${\sl . The equality is attained at the point }$x\in M$%
{\sl \ if and only if there is an orthonormal frame }$\{e_{1},...,e_{n}\}$%
{\sl \ in }$T_{x}M${\sl \ and an orthonormal frame }$\{e_{n+1},...,e_{m}\}$%
{\sl \ in }$T_{x}^{\perp }M${\sl \ in which the Weingarten operators have
the following form}

\[
A_{n+1}=\left( 
\begin{array}{lllll}
h_{11}^{n+1} & 0 & 0 & . & 0 \\ 
0 & h_{22}^{n+1} & 0 & . & 0 \\ 
0 & 0 & h_{33}^{n+1} & . & 0 \\ 
\text{.} & \text{.} & \text{.} & \text{.} & \text{.} \\ 
0 & 0 & 0 & . & h_{nn}^{n+1}
\end{array}
\right) \text{ }, 
\]
{\sl with} $h_{11}^{n+1}+h_{22}^{n+1}=h_{33}^{n+1}=...=h_{nn}^{n+1}$ {\sl and%
}

\[
A_{r}=\left( 
\begin{array}{lllll}
h_{11}^{r} & h_{12}^{r} & 0 & . & 0 \\ 
h_{12}^{r} & -h_{11}^{r} & 0 & . & 0 \\ 
0 & 0 & 0 & . & 0 \\ 
. & . & . & . & . \\ 
0 & 0 & 0 & . & 0
\end{array}
\right) \text{ },\text{ }r\in \overline{n+2,m}. 
\]

C{\footnotesize OROLLARY} 2.1.{\bf \ }{\sl If the Riemannian manifold }$%
(M,g) ${\sl , of dimension }$n\geq 3,${\sl \ admit a minimal isometric
immersion into a real space form }$\widetilde{M}(c),${\sl \ then}

\[
k\geq \tau -\frac{(n-2)(n+1)c}{2}. 
\]

\medskip By using the constrained extremum method we give another Riemannian
obstruction to minimal isometric immersions in real space forms

\[
k\leq -\tau +\frac{(n^{2}-n+2)c}{2}. 
\]

\begin{center}
\bigskip

{\bf 3. A NEW OBSTRUCTION TO MINIMAL ISOMETRIC IMMERSIONS INTO A REAL SPACE
FORM}
\end{center}

{\bf \medskip }

\quad ~Let $(M,g)$ be a Riemannian manifold of dimension $n.$ We define the
following invariants

\[
\delta _{M}^{a}=\left\{ 
\begin{array}{c}
\tau -a\min k\text{, for }0\leq a<1, \\ 
\tau -a\max k\text{, for }-1<a<0
\end{array}
,\right. 
\]

\[
\delta _{M}^{^{\prime }}=\tau +\max k,\newline
\]
where $\tau $ is the scalare curvature, and $k$ is the sectional curvature.

\medskip \medskip \medskip

With these ingredients we obtain

\medskip

\medskip

T{\footnotesize HEOREM} 3.1.{\sl \ For any real number }$a\in (-1,1),${\sl \
the invariant }$\delta _{M}^{a}${\sl \ of a Riemannian submanifold }$(M,g)$%
{\sl , of dimension }$n\geq 3,${\sl \ into a real space form }$\widetilde{M}%
(c)${\sl , of dimension }$m,${\sl \ verifies the inequality} 
\[
\delta _{M}^{a}\leq \frac{(n^{2}-n-2a)c}{2}+\frac{n(a+1)-3a-1}{n(a+1)-2a}%
\frac{n^{2}\left\| H\right\| ^{2}}{2}, 
\]
{\sl where }$H${\sl \ is the mean curvature vector of submanifold }$M${\sl \
in }$\widetilde{M}(c)${\sl . The equality is attained at the point }$x\in M$%
{\sl \ if and only if there is an orthonormal frame }$\{e_{1},...,e_{n}\}$%
{\sl \ in }$T_{x}M${\sl \ and an orthonormal frame }$\{e_{n+1},...,e_{m}\}$%
{\sl \ in }$T_{x}^{\perp }M${\sl \ in which the Weingarten operators have
the following form}

$\ $

\[
A_{r}=\left( 
\begin{array}{ccccc}
h_{11}^{r} & 0 & 0 & . & 0 \\ 
0 & h_{22}^{r} & 0 & . & 0 \\ 
0 & 0 & h_{33}^{r} & . & 0 \\ 
. & . & . & . & . \\ 
0 & 0 & 0 & . & h_{nn}^{r}
\end{array}
\right) , 
\]
{\sl with }$(a+1)h_{11}^{r}=(a+1)h_{22}^{r}=h_{33}^{r}=...=h_{nn}^{r},$ $%
\forall $ $r\in \overline{n+1,m}.$

\medskip

{\sl Proof}{\bf . }Consider $x\in M$, $\{e_{1},e_{2},...,e_{n}\}$ an
orthonormal frame in $T_{x}M$ , $\{e_{n+1},e_{n+2},...,e_{m}\}$ an
orthonormal frame in $T_{x}^{\bot }M$ and $a$ $\in (-1,1).$

\medskip

From Gauss equation it follows

$\tau -ak(e_{1}\wedge e_{2})=\frac{(n^{2}-n-2a)c}{2}+\sum\limits_{r=n+1}^{m}%
\sum\limits_{1\leq i<j\leq n}(h_{ii}^{r}h_{jj}^{r}-(h_{ij}^{r})^{2})-$

\[
-a\sum\limits_{r=n+1}^{m}(h_{11}^{r}h_{22}^{r}-(h_{12}^{r})^{2}). 
\]

\medskip Using the fact that $a$ $\in (-1,1)$, we obtain\newline
(1) $\tau -ak(e_{1}\wedge e_{2})\leq \frac{(n^{2}-n-2a)c}{2}%
+\sum\limits_{r=n+1}^{m}\sum\limits_{1\leq i<j\leq
n}h_{ii}^{r}h_{jj}^{r}-a\sum\limits_{r=n+1}^{m}h_{11}^{r}h_{22}^{r}.$

\medskip

For $r\in \overline{n+1,m}$, let us consider the quadratic form

\ 

\[
\begin{array}{c}
f_{r}:R^{n}\rightarrow R, \\ 
f_{r}(h_{11}^{r},h_{22}^{r},...,h_{nn}^{r})=\sum\limits_{1\leq i<j\leq
n}(h_{ii}^{r}h_{jj}^{r})-ah_{11}^{r}h_{22}^{r}
\end{array}
\]
and the constrained extremum problem

\[
\max f_{r}\text{ ,} 
\]
\[
\text{subject to }P:h_{11}^{r}+h_{22}^{r}+...+h_{nn}^{r}=k^{r}, 
\]
where $k^{r}$ is a real constant.

The first three partial derivatives of the function $f_{r}$ are\newline
(2) $\frac{\partial f_{r}}{\partial h_{11}^{r}}=\sum\limits_{2\leq j\leq
n}h_{jj\text{ }}^{r}-ah_{22}^{r},$\newline
$(3)$ $\frac{\partial f_{r}}{\partial h_{22}^{r}}=\sum\limits_{j\in 
\overline{1,n}\backslash \{2\}}h_{jj\text{ }}^{r}-ah_{11}^{r},$\newline
$(4)$ $\frac{\partial f_{r}}{\partial h_{33}^{r}}=\sum\limits_{j\in 
\overline{1,n}\backslash \{3\}}h_{jj\text{ }}^{r}.$

\medskip \medskip

As for a solution $(h_{11}^{r},h_{22}^{r},...,h_{nn}^{r})$ of the problem in
question, the vector $($grad) $(f_{1})$ is normal at $P$, from (2) and (3)
we obtain $\dsum\limits_{j=1}^{n}h_{jj\text{ }}^{r}-h_{11}^{r}-$\newline
$-ah_{22}^{r}=\dsum\limits_{j=1}^{n}h_{jj\text{ }%
}^{r}-h_{22}^{r}-ah_{11}^{r} $, therefore\newline
(5) $h_{11}^{r}=h_{22}^{r}=b^{r}.$

$\medskip $

From (2) and (4), it follows $\dsum\limits_{j=1}^{n}h_{jj\text{ }%
}^{r}-h_{11}^{r}-ah_{22}^{r}=\dsum\limits_{j=1}^{n}h_{jj\text{ }%
}^{r}-h_{33}^{r}.$ By using (5) we obtain $h_{33}^{r}=b^{r}(a+1)$. Similarly
one gets\newline
(6) $h_{jj}^{r}=b^{r}(a+1),$ $\forall $ $j\in \overline{3,n}.$

\medskip

As $h_{11}^{r}+h_{22}^{r}+...+h_{nn}^{r}=k^{r},$ from (5) and (6) we obtain%
\newline
(7) $b^{r}=\frac{k^{r}}{n(a+1)-2a}$\ .

\medskip We fix an arbitrary point $p\in P.$

The 2-form $\alpha $ : $T_{p}P\times T_{p}P\rightarrow R$ has the expression

\[
\alpha (X,Y)=\text{Hess}_{f_{r}}(X,Y)+\left\langle h^{\prime }(X,Y)\text{,}(%
\text{grad f}_{\text{r}})(p)\right\rangle , 
\]
where $h^{\prime }$ is the second fundamental form of $P$ in $R^{n}$ and $%
\left\langle \text{ },\text{ }\right\rangle $ is the standard inner-product
on $R^{n}$.

\ 

In the standard frame of $R^{n},$ the Hessian of $f_{r}$ has the matrix

\ 

\[
\text{Hess}_{f_{r}}=\left( 
\begin{array}{lllll}
0 & 1-a & 1 & . & 1 \\ 
1-a & 0 & 1 & . & 1 \\ 
1 & 1 & 0 & . & 1 \\ 
. & . & . & . & . \\ 
1 & 1 & 1 & . & 0
\end{array}
\right) . 
\]

\ As $P$ is totally geodesic in $R^{n}$ , considering a vector $X$ tangent
to $P$ at the arbitrary point $p$, that is, verifying the relation $%
\sum\limits_{i=1}^{n}X^{i}=0$, we have $\alpha (X,X)=2\sum\limits_{1\leq
i<j\leq
n}X^{i}X^{j}-2aX^{1}X^{2}=(\sum\limits_{i=1}^{n}X^{i})^{2}-\sum%
\limits_{i=1}^{n}(X^{i})^{2}-2aX^{1}X^{2}=$\newline
$=-\sum%
\limits_{i=1}^{n}(X^{i})^{2}-a(X^{1}+X^{2})^{2}+a(X^{1})^{2}+a(X^{2})^{2}=$%
\newline
$=-\sum%
\limits_{i=3}^{n}(X^{i})^{2}-a(X^{1}+X^{2})^{2}-(1-a)(X^{1})^{2}-(1-a)(X^{2})^{2}\leq 0. 
$ Therefore the point $(h_{11}^{r},h_{22}^{r},...,h_{nn}^{r})$, which
satisfies (5), (6), (7) is a maximum point.

\medskip

From (5) and (6) it follows\newline
(8) $f_{r}\leq
(b^{r})^{2}+2b^{r}(n-2)b^{r}(a+1)+C_{n-2}^{2}(b^{r})^{2}(a+1)^{2}-a(b^{r})^{2}= 
$

\[
=\frac{(b^{r})^{2}}{2}[n^{2}(a+1)^{2}-n(a+1)(5a+1)+6a^{2}+2a]= 
\]
\[
=\frac{(b^{r})^{2}}{2}[n(a+1)-3a-1][n(a+1)-2a]. 
\]

\medskip

By using (7) and (8), we obtain\newline
(9) $f_{r}\leq \frac{(k^{r})^{2}}{2[n(a+1)-2a]}[n(a+1)-3a-1]=\frac{%
n^{2}(H^{r})^{2}}{2}\frac{n(a+1)-3a-1}{n(a+1)-2a}.$

\medskip

The relations (1) and (9) imply\newline
(10) $\tau -ak(e_{1}\wedge e_{2})\leq \frac{(n^{2}-n-2a)c}{2}+\frac{%
n(a+1)-3a-1}{n(a+1)-2a}\frac{n^{2}\left\| H\right\| ^{2}}{2}.$

\medskip

In (10) we have equality if and only if the same thing occurs in the
inequality (1) and, in addition, (5) and (6) occurs. Therefore\newline
(11) $h_{ij}^{r}$ $=$ $0$, $\forall $ $r\in \overline{n+1,m}$, $\forall $ $%
i,j\in \overline{1,n},$ with $i\neq j$ and\newline
(12) $(a+1)h_{11}^{r}=(a+1)h_{22}^{r}=h_{33}^{r}=...=h_{nn}^{r},$ $\forall $ 
$r\in \overline{n+1,m}$.

\medskip

The relations (10), (11) and (12) imply the conclusion of the theorem.

\medskip

{\sl Remark}{\bf .} i) Making $a$ to converge at $1$ in previous inequality,
we obtain {\sl Chen's Inequality.} The conditions for which we have equality
are obtained in [1] and [6].

\medskip

ii) For $a=0$ we obtain the well-known inequality 
\[
\tau \leq \frac{n(n-1)}{2}(\left\| H\right\| ^{2}+c). 
\]

The equality is attained at the point{\sl \ }$x\in M${\sl \ }if and only if $%
x$ is a totally umbilical point.

\medskip

iii) Making $a$ to converge at $-1$ in previous inequality, we obtain 
\[
\delta _{M}^{^{\prime }}\leq \frac{(n^{2}-n+2)c}{2}+\frac{n^{2}\left\|
H\right\| ^{2}}{2}. 
\]

\medskip

T{\footnotesize HEOREM} 3{\bf .}2. {\sl The invariant }$\delta
_{M}^{^{\prime }}${\sl \ of a Riemannian submanifold }$(M,g)${\sl , of
dimension }$n\geq 3,${\sl \ into a real space form }$\widetilde{M}(c)${\sl ,
of dimension }$m,${\sl \ verifies the inequality} 
\[
\delta _{M}^{^{\prime }}\leq \frac{(n^{2}-n+2)c}{2}+\frac{n^{2}\left\|
H\right\| ^{2}}{2}, 
\]
{\sl where }$H${\sl \ is the mean curvature vector of submanifold }$M${\sl \
in }$\widetilde{M}(c)${\sl . The equality is attained at the point }$x\in M$%
{\sl \ if and only if there is an orthonormal frame }$\{e_{1},...,e_{n}\}$%
{\sl \ in }$T_{x}M${\sl \ and an orthonormal frame }$\{e_{n+1},...,e_{m}\}$%
{\sl \ in }$T_{x}^{\perp }M${\sl \ in which the Weingarten operators have
the following form}

$\ $

\[
A_{r}=\left( 
\begin{array}{ccccc}
h_{11}^{r} & 0 & 0 & . & 0 \\ 
0 & h_{22}^{r} & 0 & . & 0 \\ 
0 & 0 & 0 & . & 0 \\ 
. & . & . & . & . \\ 
0 & 0 & 0 & . & 0
\end{array}
\right) , 
\]
{\sl with }$h_{11}^{r}=h_{22}^{r},$ $\forall $ $r\in \overline{n+1,m}.$

\medskip

{\sl Proof}{\bf . }We consider the point{\bf \ }$x\in M$, the orthonormal
frames $\{e_{1},...,e_{n}\}$ in $T_{x}M$ and $\{e_{n+1},...,e_{m}\}$ in $%
T_{x}^{\perp }M$, $\{e_{1},e_{2}\}$ being an orthonormal frame in the $2-$
plane which maximize the sectional curvature at the point $x$ in $T_{x}M.$

The invariant $\delta _{M}^{^{\prime }}$ verifies\newline
(1) $\delta _{M}^{^{\prime }}=\frac{(n^{2}-n+2)c}{2}+\sum\limits_{r=n+1}^{m}%
\sum\limits_{1\leq i<j\leq
n}(h_{ii}^{r}h_{jj}^{r}-(h_{ij}^{r})^{2})+\sum%
\limits_{r=n+1}^{m}(h_{11}^{r}h_{22}^{r}-(h_{12}^{r})^{2})\leq $\newline
$\leq \frac{(n^{2}-n+2)c}{2}+\sum\limits_{r=n+1}^{m}\sum\limits_{1\leq
i<j\leq n}(h_{ii}^{r}h_{jj}^{r})+h_{11}^{r}h_{22}^{r}.$

\medskip \medskip

For $r\in \overline{n+1,m}$, let us consider the quadratic form

\ 

\[
\begin{array}{c}
f_{r}:R^{n}\rightarrow R, \\ 
f_{r}(h_{11}^{r},h_{22}^{r},...,h_{nn}^{r})=\sum\limits_{1\leq i<j\leq
n}(h_{ii}^{r}h_{jj}^{r})+h_{11}^{r}h_{22}^{r}
\end{array}
\]
and the constrained extremum problem

\[
\max f_{r}\text{ ,} 
\]
\[
\text{subject to }P:h_{11}^{r}+h_{22}^{r}+...+h_{nn}^{r}=k^{r}, 
\]
where $k^{r}$ is a real constant.

The first three partial derivatives of the function $f_{r}$ are\newline
(2) $\frac{\partial f_{r}}{\partial h_{11}^{r}}=\sum\limits_{2\leq j\leq
n}h_{jj\text{ }}^{r}+h_{22}^{r},$\newline
$(3)$ $\frac{\partial f_{r}}{\partial h_{22}^{r}}=\sum\limits_{j\in 
\overline{1,n}\backslash \{2\}}h_{jj\text{ }}^{r}+h_{11}^{r},$\newline
$(4)$ $\frac{\partial f_{r}}{\partial h_{33}^{r}}=\sum\limits_{j\in 
\overline{1,n}\backslash \{3\}}h_{jj\text{ }}^{r}.$

\medskip \medskip

As for a solution $(h_{11}^{r},h_{22}^{r},...,h_{nn}^{r})$ of the problem in
question, the vector $($grad) $(f_{1})$ is normal at $P$, from (2) and (3)
we obtain $\dsum\limits_{j=1}^{n}h_{jj\text{ }}^{r}-h_{11}^{r}+$\newline
$+h_{22}^{r}=\dsum\limits_{j=1}^{n}h_{jj\text{ }}^{r}-h_{22}^{r}+h_{11}^{r}$%
, therefore\newline
(5) $h_{11}^{r}=h_{22}^{r}=b^{r}.$

$\medskip $

From (2) and (4), it follows $\dsum\limits_{j=1}^{n}h_{jj\text{ }%
}^{r}-h_{11}^{r}+h_{22}^{r}=\dsum\limits_{j=1}^{n}h_{jj\text{ }%
}^{r}-h_{33}^{r}.$ By using (5) we obtain $h_{33}^{r}=0$. Similarly one gets%
\newline
(6) $h_{jj}^{r}=0,$ $\forall $ $j\in \overline{3,n}.$

\medskip

As $h_{11}^{r}+h_{22}^{r}+...+h_{nn}^{r}=k^{r},$ from (5) and (6) we obtain%
\newline
(7) $b^{r}=\frac{k^{r}}{2}$\ .

\medskip We fix an arbitrary point $p\in P.$

The 2-form $\alpha $ : $T_{p}P\times T_{p}P\rightarrow R$ has the expression

\[
\alpha (X,Y)=\text{Hess}_{f_{r}}(X,Y)+\left\langle h^{\prime }(X,Y)\text{,}(%
\text{grad f}_{\text{r}})(p)\right\rangle , 
\]
where $h^{\prime }$ is the second fundamental form of $P$ in $R^{n}$ and $%
\left\langle \text{ },\text{ }\right\rangle $ is the standard inner-product
on $R^{n}$.

\ 

In the standard frame of $R^{n},$ the Hessian of $f_{r}$ has the matrix

\ 

\[
\text{Hess}_{f_{r}}=\left( 
\begin{array}{lllll}
0 & 2 & 1 & . & 1 \\ 
2 & 0 & 1 & . & 1 \\ 
1 & 1 & 0 & . & 1 \\ 
. & . & . & . & . \\ 
1 & 1 & 1 & . & 0
\end{array}
\right) . 
\]

\ As $P$ is totally geodesic in $R^{n}$ , considering a vector $X$ tangent
to $P$ at the arbitrary point $p$, that is, verifying the relation $%
\sum\limits_{i=1}^{n}X^{i}=0$, we have $\alpha (X,X)=2\sum\limits_{1\leq
i<j\leq
n}X^{i}X^{j}+2X^{1}X^{2}=(\sum\limits_{i=1}^{n}X^{i})^{2}-\sum%
\limits_{i=1}^{n}(X^{i})^{2}+2X^{1}X^{2}=-\sum%
\limits_{i=3}^{n}(X^{i})^{2}-(X^{1}-X^{2})^{2}.$

\medskip

The 2-form $\alpha $ is semipositive definite. Therefore the point $%
(h_{11}^{r},h_{22}^{r},...,h_{nn}^{r})$ which satisfies (5), (6) and (7) is
a global maximum point. Using this fact and (5), (6) and (7) we obtain 
\[
\delta _{M}^{^{\prime }}\leq \frac{(n^{2}-n+2)c}{2}+\frac{n^{2}\left\|
H\right\| ^{2}}{2} 
\]

\medskip

The relation 
\[
\delta _{M}^{^{\prime }}=\frac{(n^{2}-n+2)c}{2}+\frac{n^{2}\left\| H\right\|
^{2}}{2} 
\]
occurs if and only if we have\newline
(9) $h_{ij}^{r}$ $=$ $0$, $\forall $ $r\in \overline{n+1,m}$, $\forall $ $%
i,j\in \overline{1,n}$ with $i\neq j$,\newline
(10) $h_{11}^{r}=h_{22}^{r},$ $\forall $ $r\in \overline{n+1,m}$,\newline
(11) $h_{33}^{r}=...=h_{nn}^{r}=0,$ $\forall $ $r\in \overline{n+1,m}$.

Therefore, there is an orthonormal frame{\sl \ }$\{e_{1},...,e_{n}\}${\sl \ }%
in $T_{x}M$ and an orthonormal frame $\{e_{n+1},...,e_{m}\}$ in $%
T_{x}^{\perp }M$ in which the Weingarten operators have the following form

$\ $

\[
A_{r}=\left( 
\begin{array}{ccccc}
h_{11}^{r} & 0 & 0 & . & 0 \\ 
0 & h_{22}^{r} & 0 & . & 0 \\ 
0 & 0 & 0 & . & 0 \\ 
. & . & . & . & . \\ 
0 & 0 & 0 & . & 0
\end{array}
\right) , 
\]
with{\sl \ }$h_{11}^{r}=h_{22}^{r},$ $\forall $ $r\in \overline{n+1,m}.$

\medskip

C{\footnotesize OROLLARY}{\bf \ }3.1. {\sl If the Riemannian manifold }$%
(M,g) ${\sl , of dimension }$n\geq 3,${\sl \ admit a minimal isometric
immersion into a real space form }$\widetilde{M}(c),${\sl \ then}

\[
\tau -\frac{(n-2)(n+1)c}{2}\leq k\leq -\tau +\frac{(n^{2}-n+2)c}{2}. 
\]

\medskip

C{\footnotesize OROLLARY} 3.2. {\sl If the Riemannian manifold }$(M,g)${\sl %
, of dimension }$n\geq 3,${\sl \ admit a minimal isometric immersion into a
Euclidean space}$,${\sl \ then}

\[
\tau \leq k\leq -\tau . 
\]

\medskip

\end{document}